\documentclass[12pt]{amsart}

\newcommand{\me}{\author[Quigg]{John Quigg}
\address{Department of Mathematics and Statistics
\\Arizona State University
\\Tempe, Arizona 85287}
\email{quigg@asu.edu}}

\newcommand{\kaz}{\author[Kaliszewski]{S. Kaliszewski}
\address{Department of Mathematics and Statistics
\\Arizona State University
\\Tempe, Arizona 85287}
\email{kaliszewski@asu.edu}}

\newcommand{\magnus}{\author[Landstad]{Magnus~B. Landstad}
\address{Department of Mathematical Sciences\\
Norwegian University of Science and Technology\\
NO-7491 Trondheim, Norway}
\email{magnusla@math.ntnu.no}}

\usepackage{verbatim}
\usepackage{amssymb}
\usepackage{upref}
\usepackage[all]{xy}
\usepackage[lite]{amsrefs}

\newtheorem{thm}{Theorem}[section]
\newtheorem{lem}[thm]{Lemma}
\newtheorem{cor}[thm]{Corollary}
\newtheorem{prop}[thm]{Proposition}

\theoremstyle{definition}

\newtheorem{q}[thm]{Question}

\newtheorem{ex}[thm]{Example}

\newtheorem*{notnonly}{Notation}

\newtheorem{rem}[thm]{Remark}

\numberwithin{equation}{section}

\newcommand{\secref}[1]{Section~\textup{\ref{#1}}}
\newcommand{\subsecref}[1]{Subsection~\textup{\ref{#1}}}
\newcommand{\thmref}[1]{Theorem~\textup{\ref{#1}}}
\newcommand{\corref}[1]{Corollary~\textup{\ref{#1}}}
\newcommand{\lemref}[1]{Lemma~\textup{\ref{#1}}}
\newcommand{\propref}[1]{Proposition~\textup{\ref{#1}}}

\newcommand{\exref}[1]{Example~\textup{\ref{#1}}}

\newcommand{\righttext}[1]{\qquad\text{#1 }}
\renewcommand{\for}{\righttext{for}}

\newcommand{\midtext}[1]{\quad\text{ #1 }\quad}

\renewcommand{\and}{\midtext{and}}

\newcommand{\ie}{\emph{i.e.}}
\newcommand{\eg}{\emph{e.g.}}
\newcommand{\etc}{\emph{etc.}}

\newcommand{\N}{\mathbb N}
\newcommand{\Z}{\mathbb Z}
\newcommand{\Q}{\mathbb Q}

\newcommand{\C}{\mathbb C}

\newcommand{\HH}{\mathcal H}
\newcommand{\GG}{\mathcal G}

\newcommand{\Gb}{\bar G}
\newcommand{\Hb}{\bar H}
\newcommand{\Nb}{\bar N}
\newcommand{\Mb}{\bar M}
\newcommand{\Rb}{\bar R}
\newcommand{\Qb}{\bar Q}

\newcommand{\M}{\mathrm{M}}

\newcommand{\SL}{\mathrm{SL}}
\newcommand{\GL}{\mathrm{GL}}
\newcommand{\GLtwoZ}{\mathrm{GL}(2,\Z)} 

\renewcommand{\a}{\alpha}
\renewcommand{\b}{\beta}

\newcommand{\p}{\phi}
\renewcommand{\epsilon}{\varepsilon}

\newcommand{\x}{\xi}

\newcommand{\s}{\sigma}
\newcommand{\g}{\gamma}
\renewcommand{\t}{\theta}
\newcommand{\Chi}{\raisebox{2pt}{\ensuremath{\chi}}}
\renewcommand{\epsilon}{\varepsilon}

\renewcommand{\P}{\Phi}
\renewcommand{\O}{\Omega}

\DeclareMathOperator{\aut}{Aut}
\DeclareMathOperator{\ad}{Ad}

\DeclareMathOperator*{\spn}{span}
\DeclareMathOperator*{\clspn}{\overline{\spn}}
\DeclareMathOperator*{\invlim}{\varprojlim}

\renewcommand{\>}{\rangle}

\newcommand{\inv}{^{-1}}
\newcommand{\ideal}{\vartriangleleft}

\renewcommand{\bar}{\overline}
\newcommand{\what}{\widehat}
\newcommand{\wilde}{\widetilde}

\newcommand{\cc}{\mathcal}

\renewcommand{\GG}{\bar G}
\renewcommand{\HH}{\bar H}
\newcommand{\NN}{\bar N}
\newcommand{\MM}{\bar M}
\newcommand{\RR}{\bar R}
\newcommand{\QQ}{\bar Q}

\newcommand{\F}{\cc F}
\newcommand{\E}{\cc E}

\renewcommand{\H}{\cc H}

\begin{document}

\title[Hecke $C^*$-algebras and semidirect products]
{Hecke $C^*$-algebras and semidirect products}

\kaz
\magnus
\me

\subjclass[2000]{Primary 46L55; Secondary 20C08}

\keywords{Hecke algebra,
group $C^*$-algebra,
Morita equivalence,
semidirect product}

\begin{abstract}
We analyze Hecke pairs $(G,H)$ and the associated Hecke algebra
$\cc H$
when $G$ is a semidirect product $N\rtimes Q$ and $H=M\rtimes R$
for subgroups $M\subset N$ and $R\subset Q$ with $M$ normal in~$N$.
Our main result shows that 
when $(G,H)$ coincides with its Schlichting completion
and $R$ is normal in~$Q$, 
the closure of~$\cc H$ in~$C^*(G)$ is Morita-Rieffel
equivalent
to a crossed product $I\rtimes_\beta Q/R$, where~$I$ is a certain ideal
in the fixed-point algebra $C^*(N)^R$.  Several concrete examples are
given
illustrating and applying our techniques, including some involving
subgroups of $\GL(2,K)$ acting on~$K^2$, where $K=\Q$ or
$K=\Z[p^{-1}]$.
In particular we look  at the $ax+b$-group of a quadratic
extension of~$K$.
\end{abstract}


\maketitle

\section*{Introduction}
\label{intro}

A \emph{Hecke pair} $(G,H)$ comprises a group $G$ and a
subgroup
$H$ for which every double coset is a finite union of left cosets, and
the
associated \emph{Hecke algebra}, generated by the characteristic
functions
of double cosets, reduces to the group $*$-algebra
of $G/H$
when $H$ is normal.

In \cite{hecke} we introduced the \emph{Schlichting completion} 
$(\Gb,\Hb)$ of the Hecke pair $(G,H)$
as a
tool for analyzing Hecke algebras, based in part upon work of Tzanev
\cite{tza}. 
(A slight variation on this construction
appears in \cite{willis}.)
The idea is that $\HH$ is a compact open subgroup of $\GG$ such
that the Hecke algebra of $(\GG,\HH)$ is naturally identified
with the Hecke algebra $\H$ of $(G,H)$.  The characteristic function
$p$ of $\HH$ is a projection in
the group $C^*$-algebra
$A:=C^*(\GG)$, and 
$\H$ can be identified with $pC_c(\GG)p\subset A$;
thus the closure of $\H$ in $A$
coincides with the corner $pAp$, which 
is Morita-Rieffel equivalent to the ideal $\bar{ApA}$.

In \cite{hecke} we were mainly interested in studying when
$pAp$ is the
enveloping $C^*$-algebra of the Hecke algebra $\H$, and when
the projection $p$ is full in $A$,
making the $C^*$-completion $pAp$ of $\H$ Morita-Rieffel
equivalent to the group $C^*$-algebra $A$.
We had the most success when $G=N\rtimes Q$ was a semidirect product with
the Hecke subgroup~$H$ contained in the normal subgroup $N$. 

In this paper we again
consider $G=N\rtimes Q$, but now we allow 
$H=M\rtimes R$, where $M$ is a normal subgroup of $N$
and $R$ is a subgroup of $Q$ which normalizes $M$.  
Briefly:
\begin{equation}\label{setup}
\begin{matrix}
G & = & N & \rtimes & Q\\
\vee & & \triangledown & & \vee\\
H & = & M & \rtimes & R.
\end{matrix}
\end{equation}
This leads to a refinement of the
Morita-Rieffel equivalence $\bar{ApA}\sim pAp$ (see \thmref{cross
thm}).

We begin in \secref{prelim} by recalling our conventions from
\cite{hecke} regarding Hecke algebras.  
In \secref{group} we describe the main properties of
our group-theoretic setup~\eqref{setup}.
In particular, we characterize the reduced Hecke pairs
in terms of $N$, $Q$, $M$, and $R$.

In order to effectively analyze how our semidirect-product
decomposition affects the Hecke topology, we need to go into somewhat
more detail than might be expected.
In particular, we must exercise some care
to obtain the
semidirect-product
decomposition 
\[\Gb=\Nb\rtimes \Qb,\quad
\Hb=\Mb\rtimes \Rb\]
for the Schlichting completion (see \corref{semidirect
completion}), and to 
describe various bits of this completion
as inverse limits of groups (see \thmref{inverse}).

\secref{crossed} is preparatory for \secref{hecke crossed}, but the
results may be of
independent interest.
In \propref{newrosenberg} we show that if
$(B,Q,\a)$ is an action, $R$ is a compact normal subgroup of $Q$
and $(B^R,Q/R,\b)$ is the associated action,
then the projection
$q=\int_R r\,dr$ is in $M(B\times_\a Q)$ and
$B^R\times_\b Q/R\cong q(B\times_\a Q)q$. This generalises the
result of
\cite{rosenberg}.

We also show in \thmref{combine} that under this correspondence the
ideal
$ \bar{(B\times_\a G)p(B\times_\a G)}$ is mapped to an ideal
$I\times_\b Q/R $ where $I$ is a $Q/R$-invariant ideal of $B^R$.

In \secref{hecke crossed},
we assume that $R$ is normal in $Q$,
and (without loss of generality)
that the pair $(G,H)$ is equal to
its Schlichting completion.
The main result is \thmref{cross thm}, in which we take full advantage
of the
semidirect-product decomposition to show that
the Hecke $C^*$-algebra $p_HC^*(G)p_H$ is Morita-Rieffel equivalent to a
crossed
product
$I\times_\b Q/R $, where $I$ is the ideal in $C^*(N)^R $ generated by
$\{\a_s(p_M):s\in Q\}$.
We look briefly at the special case where the normal subgroup $N$ is
abelian.

Finally, in \secref{example} we give some examples to illustrate our
results.
Classical Hecke algebras have most commonly treated pairs of
semi-simple groups such
as $(\GL(n,\Q),\SL(n,\Z))$. The work of Bost and Connes \cite{bc} showed
the
importance of also studying Hecke pairs of solvable groups. In the
examples we mostly deal with the following situation:
$K$ is either the field $\Q$ of rational numbers or
the field $\Z[p\inv]$ of rational numbers  with denominators of the
form $p^n$; $N=K^2$; $M=\Z^2$;
$Q$ is a subgroup of $\GL (2,K)$ containing the diagonal subgroup,
acting on $N$ in the obvious way;
and $R=Q\cap \GLtwoZ$.
It is not so difficult to see
that the Schlichting completions are $p$-adic or adelic versions of the
same groups.

As to specific examples we look at the algebra studied by 
Connes-Marcolli
in \cite{connes-marcolli}, see also \cite{lln:hecke}. Here $R$ is not
normal in $Q$, so the full results of \secref{hecke crossed} do not
apply.
On the other hand,
if $R$ is normal in $Q$ then \corref{M directed} does apply, and as in
\cite{lr:ideal} one can use the Mackey orbit method to study the ideal
structure of the $C^*$-algebras involved.
A particular example of this is the $ax+b$-group over a quadratic
extension $K[\sqrt d]$ treated in \cite {laca-franken}, and we shall see
that this example raises some interesting questions. We also look at a
nilpotent example, \ie, one version of the Heisenberg group over the
rationals.

After we had completed the research for this paper, we became aware of
the recent preprint \cite{lln:hecke}, which treats semidirect-product
Hecke pairs in a way quite similar to ours. 
The present paper and \cite{lln:hecke} 
were written independently, and the techniques
have only incidental overlap.
We should mention that we treat only the
case where $M$ is normal in $N$, while the context in \cite{lln:hecke}
seems to be more general. 
Thus, for example, it would be difficult to adapt our results on
inverse limits (see \subsecref{inverse limit}) to the context of
\cite{lln:hecke}.

We would like to thank Arizona State University, the
Norwegian University of Science and Technology, and the Norwegian
Science Foundation, who
all have supported this research.
We are also grateful to the referee for suggesting many improvements.

\section{Preliminaries}
\label{prelim}

We adopt the conventions of \cite{hecke}, which contains more
references.
A \emph{Hecke pair} $(G,H)$ comprises a group $G$ and a 
\emph{Hecke subgroup}
$H$, \ie, one for which every double coset $HxH$ is a finite union of
left cosets $\{y_1H,\dots,y_{L(x)}H\}$.
A good reference for the basic theory of Hecke pairs is \cite{kri}.
A Hecke pair $(G,H)$ is \emph{reduced} if $\bigcap_{x\in
G}xHx\inv=\{e\}$, and a reduced Hecke pair $(G,H)$ is a
\emph{Schlichting pair} if $G$ is locally compact Hausdorff and $H$
is compact and open in $G$.
In \cite{hecke}*{Theorem~3.8}, we gave a new proof of
\cite{tza}*{Proposition~4.1}, which says that every reduced Hecke
pair $(G,H)$ can be embedded in an essentially unique Schlichting pair
$(\GG,\HH)$, which we call the \emph{Schlichting completion} of
$(G,H)$.
Specifically, $\GG$ is the completion of $G$ in the
(two-sided uniformity defined by the) \emph{Hecke topology} having a
local subbase $\{xHx\inv\mid x\in G\}$ of neighborhoods of $e$,
and $(\GG,\HH)$ is unique in the sense that if $(L,K)$ is any
Schlichting pair and $\s:G\to L$ is a homomorphism such that $\s(G)$
is dense and $H=\s\inv(K)$, then $\s$ extends uniquely to a
topological isomorphism $\bar\s:\GG\to L$, and moreover
$\bar\s(\HH)=K$.

The associated \emph{Hecke algebra} is the vector subspace $\H$ of
$\C^G$ spanned by the characteristic functions of double $H$-cosets,
with operations defined by
\begin{align*}
f*g(x)
&=\sum_{yH\in G/H}f(y)g(y\inv x)
\\
f^*(x)
&=\bar{f(x\inv)}\Delta(x\inv),
\end{align*}
where $\Delta(x)=L(x)/L(x^{-1})$
and $L(x)$ is the number of left cosets $yH$ in the double coset
$HxH$.
Warning: some authors do not include the factor of $\Delta$ in the
involution; for us it arises naturally when we embed $\H$ in
$C_c(\GG)$ 
(see \cite{hecke}*{Section~1}). 
One way to see how this
embedding goes is the
following: let $p=\Chi_{\HH}$, which is a projection in
$C_c(\GG)$ when the Haar measure on $\GG$ has been normalized so
that $\HH$ has measure~$1$. Then
the restriction map $f\mapsto f|G$ gives a $*$-isomorphism of the
convolution algebra $pC_c(\GG)p$ onto $\H$.

\subsection*{Notation}
$H<G$ means $H$ is a subgroup of $G$.
$H\ideal G$ means $H$ is a normal subgroup of $G$.
If $N\ideal G$ and $Q<G$ such that $N\cap Q=\{e\}$
and $NQ=G$, then $G$ is the (internal) semidirect product of $N$ by
$Q$, and we write $G = N\rtimes Q$.

\section{Groups}
\label{group}

Here we describe the main properties of our group-theoretic 
setup~\eqref{setup} for Hecke semidirect products.
We need to establish many elementary facts from group theory
which are not standard, so we will give more detail than might seem
necessary.


\subsection{Generalities}

We will be interested in subgroups of $H$
of the form $LS$, where $L<M$ and $S<R$. 
Note that $LS<MR$ if and only if $S$ normalizes $L$.

\begin{lem}
If $A,B,C$ are subgroups of $G$ with:
\begin{enumerate}
\item $A\supset B$;
\item $A\cap C=\{e\}$;
\item $AC=CA$;
\item $BC=CB$,
\end{enumerate}
then
\[[AC:BC]=[A:B].\]
\end{lem}

\begin{proof}
The map $aB\mapsto aBC:A/B\to AC/BC$ is obviously well-defined and surjective, and is injective because
\[a_1BC=a_2BC
\implies a_2\inv a_1\in BC\cap A=B.\qedhere\]
\end{proof}

\begin{cor}\label{MRLS}
Suppose $L<M$ and $S<R$, and suppose $S$ normalizes~$L$, 
so that $LS$ is a subgroup of $MR$.  
Then
\[
[M:L][R:S] = [MR:LS].
\]
\end{cor}

\begin{proof}
We have
\[[MR:LS]=[MR:MS][MS:LS],\]
so the result follows from the above lemma.
\end{proof}

\begin{notnonly}
For any subgroup $K$ of $G$ and $x\in G$, we define
\[
K_x = K \cap xKx^{-1}.
\]
Thus $K_x$ is precisely the stabilizer subgroup of the coset
$xK$ under the action of $K$ on $G/K$ by \emph{left translation},
and
\begin{equation}\label{Kx}
[K:K_x] = \left| KxK / K \right|.
\end{equation}

If $T$ is another subgroup of $G$, we let
\[
T_{x,K} = \{ t\in T \mid txKt^{-1} = xK \}
\]
denote the stabilizer subgroup of $xK$ under the action of $T$ by
\emph{conjugation} on the set of all subsets of $G$; thus
\begin{equation}\label{TTxK}
[T:T_{x,K}] = \left| \{ txKt^{-1} \mid t\in T \} \right|.
\end{equation}
Note that if $T$ normalizes $K$, then 
the conjugation action of $T$ descends to $G/K$.

For $E\subset G$, we further define
\[
K_E=\bigcap_{x\in E} K_x
\midtext{and}
T_{E,K}=\bigcap_{x\in E} T_{x,K}.
\]
\end{notnonly}

It will also be useful to observe that
if $\{M_i\}_{i\in I}$ is a family of subgroups of $N$
and $\{R_i\}_{i\in I}$ is a family of subgroups of $Q$
such that $R_i$ normalizes $M_i$ for each $i\in I$, 
then, because $N\cap Q = \{e\}$, we have
\begin{equation}\label{MRMR}
\bigcap_{i\in I}M_iR_i
=\biggl(\bigcap_{i\in I}M_i\biggr)
\biggl(\bigcap_{i\in I}R_i\biggr)
\end{equation}

\begin{lem}\label{RnL}
Let $L$ be a subgroup of $N$ which is normalized by~$R$. 
For any $r\in R$ and $n\in N$, the following are equivalent\textup:
\begin{enumerate}
\item $r\in R_{n,L}$;
\item $rnr\inv\in nL$;
\item $r\in nLRn\inv$.
\end{enumerate}
\end{lem}

\begin{proof}[Sketch of Proof]
(i) $\implies$ (ii) $\implies$ (iii) is clear.
(iii) $\implies$ (ii) uses $N\cap Q=\{e\}$. 
(ii) $\implies$ (i) because $R$ normalizes $L$.
\end{proof}

Taking $L=M$ in \lemref{RnL} 
and using $H=MR$, we have
\begin{equation}\label{star}
R_{n,M} = R\cap nMRn^{-1} =  R\cap nHn^{-1}
\supset R\cap nRn^{-1} = R_n.
\end{equation}
From this we deduce:

\begin{lem}\label{HMR}
For any $n\in N$ and $q\in Q$,
\begin{enumerate}
\item $H_n = M R_{n,M}$;
\item $H_q = M_q R_q$;
\item $H_{qn} \cap H_q  = M_q ( q R_{n,M} q^{-1} \cap R)$.
\end{enumerate}
\end{lem}

\begin{proof}
(i)  Suppose $h=mr\in H_n$ for $m\in M$ and $r\in R$.  
Then
\[r\in m^{-1} nMRn^{-1} = n (n^{-1}m^{-1}nMR)n^{-1}
= nMRn^{-1} = nHn^{-1},\]
so (using \eqref{star})
\[mr \in m(R\cap nHn^{-1}) \subset MR_{n,M}.\]
Thus $H_n\subset MR_{n,M}$.
Conversely, also using \eqref{star}, 
\[MR_{n,M} = M(R\cap nHn^{-1}) \subset MR \cap MnHn^{-1}
= H \cap nHn^{-1} = H_n.\]

(ii) By \eqref{MRMR} we have
\begin{align*}
H_q 
&= H \cap qHq^{-1} 
= MR \cap (qMq^{-1})(qRq^{-1})
\\&= (M\cap qMq^{-1})(R\cap qRq^{-1}) 
= M_q R_q.
\end{align*}

(iii) Using part~(i) and \eqref{MRMR} we have
\begin{align*}
H_{qn} \cap H_q
&= H\cap qHq^{-1} \cap qnHn^{-1}q^{-1}
= H \cap q(H_n)q^{-1}
\\&= MR \cap (qMq^{-1})(qR_{n,M}q^{-1})
= M_q (R \cap qR_{n,M}q^{-1}).
\end{align*}
\end{proof}


\subsection{Hecke pairs}
Since $[H:H_x] = |HxH/H|$ for any $x\in G$, 
the pair $(G,H)$ is Hecke if and only if 
each subgroup $H_x$ has finite index in $H$.  
Applying this to the pair $(N\rtimes Q,Q)$,
we see that $(N\rtimes Q,Q)$ is Hecke if and only if
$[Q:Q_n] = [Q:Q_{n,\{e\}}]<\infty$ for each $n\in N$.
The next proposition extends this observation to our more
general context.

\begin{prop}\label{hecke}
The following are equivalent\textup:
\begin{enumerate}
\item $(G,H)$ is a Hecke pair
\item $[R:R_q]$, $[M:M_q]$ and $[R:R_{n,M}]$ are all finite for each
$q\in Q$ and $n\in N$
\item $(Q,R)$, $(G,M)$ and $(N/M\rtimes R,R)$ are Hecke pairs
\item $(Q,R)$, $(G,M)$ and $(NR,H)$ are Hecke pairs.
\end{enumerate}
\end{prop}

\begin{proof}
If $(G,H)$ is a Hecke pair, then for all $q\in Q$ and $n\in N$ we have
\[
[M:M_q][R:R_q]=[MR:M_qR_q]=[H:H_q]<\infty
\]
{and}
\[
[R:R_{n,M}]=[MR:MR_{n,M}]=[H:H_n]<\infty,
\]
so (i) implies (ii). 
Conversely, assuming~(ii), for any $q\in Q$ and $n\in N$, 
\lemref{HMR} gives
\begin{align*}
[H:H_{qn}]
&\leq [H:H_{qn}\cap H_q]
= [MR:M_q(qR_{n,M}q\inv\cap R)]\\
&=[M:M_q][R:qR_{n,M}q\inv\cap R]\\
&=[M:M_q][R:R_q][R_q:qR_{n,M}q\inv\cap R],
\end{align*}
which is finite because 
for any subgroups $S\supset T$ of $G$ we have 
$[R\cap S:R\cap T]\le [S:T]$ and $[qSq\inv:qTq\inv]=[S:T]$.
Thus~(ii) implies~(i).

If $q\in Q$ and $n\in N$ then $qnMn\inv q\inv=qMq\inv$, so
$[M:M_q]<\infty$ for all $q\in Q$
if and only if $(G,M)$ is Hecke.
As observed above, 
$R$ is a Hecke subgroup of $N/M\rtimes R$
if and only if, for each $nM\in N/M$,
the stabilizer subgroup of $nM$ in $R$ (acting by conjugation)
has finite index in $R$.
Since this subgroup is precisely $R_{n,M}$, we have
$[R:R_{n,M}]<\infty$ for all $n\in N$ 
if and only if $(N/M\rtimes R,R)$ is Hecke.
Therefore (ii) if and only if (iii).

Finally, if $n\in N$ and $r\in R$ then $nrHr\inv n\inv=nHn\inv$, so
\[
[H:H_{nr}]=[H:H_n]=[R:R_{n,M}],
\]
therefore (iii) if and only if (iv).
\end{proof}

\begin{prop}
\label{reduced}
Suppose $(G,H)$ is a Hecke pair.
Then the following are equivalent\textup:
\begin{enumerate}
\item
$(G,H)$ is reduced;
\item
$M_Q=\{e\}$ and $R_{N,\{e\}}\cap R_Q = \{e\}$.
\end{enumerate}
\end{prop}

\begin{proof}
Since $(G,H)$ is reduced if and only if $H_G=\{e\}$,
the proposition 
will follow
easily from the identity
\begin{equation}\label{HG}
H_G = M_Q (R_{N,M_Q}\cap R_Q).
\end{equation}
To establish~\eqref{HG}, we first
use \lemref{HMR} (iii)
and \corref{MRLS} to get
\begin{align*}
H_G 
&= \bigcap_{x\in G}H_x =\bigcap_{q\in Q,n\in N}H_{qn}
=\bigcap_{q\in Q,n\in N}(H_q\cap H_{qn})\\
&=\bigcap_{q\in Q,n\in N}M_q(qR_{n,M}q\inv\cap R)\\
&=\biggl(\bigcap_{q\in Q}M_q\biggr)
  \biggl(\bigcap_{q\in Q,n\in N}qR_{n,M}q\inv\cap R\biggr).
\end{align*}

Further, 
\begin{align*}
\bigcap_{q\in Q,n\in N}qR_{n,M}q\inv \cap R
&= \bigcap_{q\in Q,n\in N}R_{qnq\inv,qMq\inv}\cap R_q\\
&= \bigcap_{q\in Q,n\in N}R_{n,qMq\inv}\cap \bigcap_{q\in Q}R_q\\
&= R_{N,M_Q} \cap R_Q.
\qedhere
\end{align*}
\end{proof}

Note that
$R_{N,\{e\}}$ consists of those elements of $R$ which commute element-wise with $N$.


\subsection{Hecke topology}
In addition to our semidirect product setup~\eqref{setup}, 
now assume that $(G,H)$ is a reduced Hecke pair.
Let $(\GG,\HH)$ denote its Schlichting completion.

\begin{prop}
\label{subbase}
The relative Hecke topologies of the relevant subgroups have the
following subbases at the identity:
\begin{enumerate}
\item for both $N$ and $M$\textup: $\{M_q \mid q\in Q\}$;
\item for $Q$\textup: $\{qR_{n,M}q\inv \mid q\in Q,n\in N\}$;
\item for $R$\textup: $\{R\cap qR_{n,M}q\inv \mid q\in Q,n\in N\}$.
\end{enumerate}
\end{prop}

\begin{proof}
(i) follows from the computation
\begin{align*}
N\cap qnHn\inv q\inv
=qn(N\cap H)n\inv q\inv
=qnMn\inv q\inv
=qMq\inv
\end{align*}
and its immediate consequence, 
$M\cap qnHn\inv q\inv =M_q$.

For (ii), we have
\begin{align*}
Q\cap nHn\inv
=Q\cap nMRn\inv
\subset Q\cap MNRN
=Q\cap NR
=R,
\end{align*}
so
\[
Q\cap qnHn\inv q\inv
=q(Q\cap nHn\inv)q\inv
=q(R\cap nHn\inv)q\inv
=qR_{n,M}q\inv.
\]

Finally, (iii) follows from~(ii).
\end{proof}

The following corollary should be compared with
\cite{lln:hecke}*{Theorem~2.9(ii)};
their extra hypothesis is satisfied in
our special case ($M\vartriangleleft N$), but it would be complicated to
verify that our result follows from theirs because their construction is
significantly different from ours.

\begin{cor}
\label{semidirect completion}
If $(G,H)$ as in~\eqref{setup} is a reduced Hecke pair
with Schlichting completion $(\GG,\HH)$, then 
\[\GG=\NN\rtimes\QQ
\and \HH=\MM\rtimes\RR,\]
where the closures are all taken in $\GG$.
\end{cor}

\begin{proof}
First of all, to show that $\GG$ is the semidirect product
$\NN\rtimes\QQ$ of its subgroups $\NN$ and $\QQ$ requires:
\begin{enumerate}
\item $\NN\ideal\GG$;
\item $\GG=\NN\;\QQ$;
\item $\NN\cap\QQ=\{e\}$;
\item $\GG$ has the product topology of $\NN\times\QQ$.
\end{enumerate}
Item~(i) is obvious. To see~(ii), note that the
subgroup $\NN\;\QQ$ contains both $G=NQ$ and $\MM\;\RR$. Since $\MM$ is
compact, the subgroup $\MM\;\RR$ is closed, and it follows that
$\HH=\MM\;\RR$. This implies~(ii), since every coset in $\GG/\HH$
can be expressed in the form $x\HH$ for $x\in G$.

For~(iii), note that the quotient map $\psi:G\to
Q\subset\QQ$ is continuous for the Hecke topology of $G$ and the
relative Hecke topology of $Q$, because a typical subbasic
neighborhood of $e$ in $Q$ is of the form $qR_{n,M}q\inv$ for $q\in Q$
and $n\in N$, and
\[\psi\inv(qR_{n,M}q\inv)=NqR_{n,M}q\inv\]
contains the neighborhood
\[H_{qn}\cap H_q=M_q(R\cap qR_{n,M}q\inv)\]
of $e$ in $G$.
Since
$\QQ$ is a complete topological group, $\psi$ extends uniquely to a
continuous homomorphism $\bar\psi:\GG\to\QQ$. Because $\psi$ takes
$N$ to $e$ and agrees with the inclusion map on $Q$, by density and
continuity $\bar\psi$ takes $\NN$ to $e$ and agrees with the
inclusion map on $\QQ$. Therefore $\NN\cap\QQ=\{e\}$.

To see how~(iv)
follows, note that 
the multiplication map $(n,q)\mapsto nq$ of $\NN\times\QQ$ onto $\GG$
is continuous by definition, and
its inverse $x\mapsto
(x\bar\psi(x)\inv,\bar\psi(x))$ is also continuous because
$\bar\psi$ is, as shown above.

It only remains to show that $\HH=\MM\rtimes\RR$, but this follows
immediately: we have $\MM\cap\RR=\{e\}$, and the subgroup $\MM\;\RR$ has
the product topology since $\NN\;\QQ$ does.
\end{proof}


\subsection{Inverse limits}
\label{inverse limit}

Here we again assume that $(G,H)$ is a reduced Hecke pair.
For each of our groups $M$, $N$, $R$, $H$, and $Q$ 
we want to describe the closure as an inverse limit of groups, 
so that we capture both the algebraic and
the
topological structure. 
From \cite{hecke}*{Proposition~3.10},
we know that the closure is topologically the
inverse
limit of the coset spaces of finite intersections of stabilizer
subgroups.
To get the algebraic structure we need enough of 
these intersections to be
normal
subgroups.
In the case of $M$ and $N$, we already have what we need, since each
$M_q$ is
normal in $N$, and hence also in $M$.
However, for $R$ we need to do more work.

\begin{lem}
\label{normsbgp}
Suppose $L<M$ and $S<R$.
Then $LS\ideal MR$
if and only if
\begin{enumerate}
\item $L\ideal MR$,
\item $S\ideal R$, and
\item $S\subset R_{M,L}$.
\end{enumerate}
Moreover, in this case
\[
MR/LS\cong (M/L)\rtimes (R/S).
\]
\end{lem}

\begin{proof}
First assume $LS\ideal MR$.
Then
\[S = R\cap LS \ideal R\cap MR = R,\]
and since $M\ideal MR$, we also have
\[L = M\cap LS \ideal MR.\]
For~(iii), fix $s\in S$ and $m\in M$.  Then $m\inv sm\in LS$ because $LS\ideal MR$, so $m\inv sms\inv\in LS$. On the other hand, $m\inv sms\inv\in M$ because $S\subset R$ and $R$ normalizes $M$. Thus
\[m\inv sms\inv\in LS\cap M=L,\]
so $s\in R_{m,L}$.

Conversely, assume
(i)--(iii).
Then it suffices to show that $M$ conjugates $S$ into $LS$:
for $m\in M$ and $s\in S$ we have
$m\inv sms\inv\in L$ by \lemref{RnL}~(ii),
and hence $m\inv sm\in LS$.

For the last statement, it is routine to verify that the map
\[
mrLS\mapsto (mL,rS)
\righttext{for}m\in M,r\in R
\]
gives a well-defined isomorphism.
\end{proof}

\begin{notnonly}
For $E\subset Q$ and $F\subset N$ put
\[
R^E_F
=\bigcap_{q\in E}qR_{F,M}q\inv \cap R
=\bigcap_{q\in E}\bigcap_{n\in F}qR_{n,M}q\inv \cap R.
\]
\end{notnonly}

Note that the families
\[\{M_E:E\subset Q\text{ finite}\}
\and
\{R^E_F:\text{ both }E\subset Q\text{ and }F\subset N\text{ finite}\}\]
are neighborhood bases at~$e$ in the relative 
Hecke topology of~$M$ and~$R$, respectively.

\begin{notnonly}
Let $\E$ be the family of all subsets $E\subset Q$ such that:
\begin{enumerate}
\item $E$ is a finite union of cosets in $Q/R$;
\item $e\in E$;
\item $RE=E$,
\end{enumerate}
and let $\F$ be the family of all pairs $(E,F)$ such that:
\begin{enumerate}
\setcounter{enumi}{3}
\item $E\in \E$;
\item $F$ is a finite union of cosets in $N/M$;
\item $q\inv Mq\subset F$ for all $q\in E$.
\end{enumerate}
\end{notnonly}

\begin{lem}\label{2-10}
For all $(E,F)\in\F$:
\begin{enumerate}
\item $R^E_F\ideal R$;
\item $[R:R^E_F]<\infty$.
\end{enumerate}
\end{lem}

\begin{proof}
$R^E_F$ is a subgroup of $R$ because $R_{F,M}$ is.
For $r\in R$ we have
\[
rR^E_Fr\inv
=\bigcap_{q\in E}r(qR_{F,M}q\inv\cap R)r\inv
=\bigcap_{q\in E}rqR_{F,M}q\inv r\inv\cap R
=R^E_F
\]
since $rE=E$. This proves (i).

For (ii), first note that $[R:R_{F,M}]<\infty$ because $|F/M|<\infty$ and $R_{n,M}$ only depends upon the coset $nM$.
Thus
\[R_0:=\bigcap_{r\in R}rR_{F,M}r\inv\]
has finite index in $R$.
For each coset $tR$ contained in $E$ we have
\[
\bigcap_{q\in tR}qR_{F,M}q\inv
=\bigcap_{r\in R}trR_{F,M}r\inv t\inv
=tR_0t\inv.
\]
Thus
\[\bigcap_{q\in tR}qR_{F,M}q\inv\cap R\]
has finite index in $R$.
Letting $E=\{t_1R,\dots,t_kR\}$, it follows that
\[
\bigcap_{q\in E}qR_{F,M}q\inv\cap R
=\bigcap_{i=1}^k
\left(\bigcap_{q\in t_iR}qR_{F,M}q\inv\cap R\right)
\]
has finite index in $R$.
\end{proof}

\begin{lem}\label{2-11}
For all $E\in \E$:
\begin{enumerate}
\item $M_E\ideal N$;
\item $M_E\ideal M$;
\item $M_E\ideal H$;
\item $[M:M_E]<\infty$.
\end{enumerate}
\end{lem}

\begin{proof}
(i)~holds because $M_q\ideal N$ for each $q$,
and~(ii) follows since $M_E\subset M$. 

(iii). For $r\in R$ we have
\[rM_Er\inv
=\bigcap_{q\in E}r(qMq\inv\cap M)r\inv
=\bigcap_{q\in E}rqMq\inv r\inv\cap M
=M_E\]
since $rE=E$. Thus $M_E\ideal MR=H$ by (ii).

(iv). For each coset $tR$ contained in $E$ we have
\[\bigcap_{q\in tR}qMq\inv
=\bigcap_{r\in R}trMr\inv t\inv
=tMt\inv.\]
Thus $\bigcap_{q\in tR}M_q=M_t$ has finite index in $M$,
and it follows that 
$M_E = \bigcap_{q\in E}M_q$ has finite index in $M$ as well.
\end{proof}

\begin{lem}\label{2-12}
For all $(E,F)\in\F$ we have
\[R^E_F\subset R_{M,M_E}.\]
\end{lem}

\begin{proof}
Fix $s\in R^E_F$ and $m\in M$; we need to show that $s\in R_{m,M_E}$. Thus, for $q\in E$, we must show
\[m\inv sms\inv\in qMq\inv.\]
We have $q\inv mq\in F$, so
$s\in qR_{q\inv mq,M}q\inv$.
It follows that
\[
q\inv m\inv sms\inv q
=(q\inv m\inv q)(q\inv sq)(q\inv mq)(q\inv s\inv q)
\in M,
\]
hence $m\inv sms\inv\in qMq\inv$, as desired.
\end{proof}

Lemmas~\ref{2-10}--\ref{2-12} yield the following:

\begin{prop}
For all $(E,F)\in\F$ we have
\[M_ER^E_F\ideal H
\and
[H:M_ER^E_F]<\infty.\]
\end{prop}

\begin{thm}
\label{inverse}
With the above notation, we have:
\begin{enumerate}
\item $\bar M=\invlim_{E\in\E}M/M_E$;
\item $\bar N=\invlim_{E\in\E}N/M_E$;
\item $\bar R=\invlim_{(E,F)\in\F}R/R^E_F$;
\item $\bar H=\invlim_{(E,F)\in\F}M/M_E\rtimes R/R^E_F$,
\end{enumerate}
all as topological groups.
\end{thm}

\begin{proof}
By the preceding results, it suffices to show that for all finite subsets
\[E'\subset Q\and F'\subset N\]
there exists $(E,F)\in\F$ such that
\[M_E\subset M_{E'}\and R^E_F\subset R^{E'}_{F'}.\]
Put
\[E''=\bigl(E'\cup\{e\}\bigr)R
\and
F''=\bigl(F'\cup\{e\}\bigr)M.\]
Since $(Q,R)$ is Hecke, $E:=RE''$ is a finite union of cosets in $Q/R$, and it follows that $E\in\E$. We have
$M_E\subset M_{E'}$ {since} $E\supset E'$.

Let $M_0$ be the subgroup of $N$ generated by the conjugates $q\inv Mq$ for $q\in E$. Then $M_0\ideal N$ since $q\inv Mq\ideal N$ for each $q$. Since $E$ is a finite union of double cosets of $R$ in $Q$, and since $(Q,R)$ is Hecke, $E$ is a finite union of right cosets of $R$ in $Q$. Thus the family $\{q\inv Mq:q\in E\}$ is finite.
Since $M_0$ is the product of the subgroups $q\inv Mq$ (because they are normal in $N$),
it follows that $[M_0:M]<\infty$.
Thus, putting
$F=M_0F''$, we have
$(E,F)\in\F$, and moreover
$R^E_F\subset R^{E'}_{F'}$
{since}
$E\supset E'$
and
$F\supset F'$.
\end{proof}

As a topological space,
$\QQ=\invlim_{E,F}Q/R_F^E$, but since
the subgroups $R_F^E$ are not in general normal in $Q$, 
the group structure of $\QQ$ is more complicated. 
For details on this, we
refer to \cite{hecke}*{Remark~3.11}.
In the special case where $Q$ is abelian, 
we do have $R^E_F\ideal Q$, so
\[
\QQ=\invlim_{E,F}Q/R_F^E
\]
as topological groups.

\section{Crossed products}
\label{crossed}

In this section we prove
a few results concerning crossed products, subgroups,
and projections. We state
these results in
somewhat greater generality than we require, since they might be
useful elsewhere and no extra work is required.

\subsection*{Compact subgroups}

Let $R$ be a compact normal subgroup of a locally compact group
$Q$.
We identify $Q$ and $C_c(Q)$ with their canonical images in
$M(C^*(Q))$ and $C^*(Q)$, respectively.
Normalize the Haar measure on $R$ so that $R$ has measure $1$.
Then $q:=\Chi_R$ is a central projection in $M(C^*(Q))$,
and the map $\tau:Q/R\to M(C^*(Q))$ defined by
\begin{equation}
\label{tau}
\tau(sR)=sq\for s\in Q
\end{equation}
integrates to give an isomorphism of $C^*(Q/R)$ with the ideal
$C^*(Q)q$ of $C^*(Q)$.

Let $\a$ be an action of $Q$ on a $C^*$-algebra $B$.
We identify $B$ and $C^*(Q)$ with their canonical images in
$M(B\times_\a Q)$.
Thus $q$ is a projection in $M(B\times_\a Q)$,
and we may regard $\tau$ as a homomorphism of $Q/R$ into
$M(B\times_\a Q)$.

Let $\P(b)=\int_R\a_r(b)\,dr$ be the faithful conditional expectation
of $B$ onto the fixed-point algebra
$B^R$.
Then an elementary calculation shows that
\[
qbq=\P(b)q=q\P(b)\for b\in B
\]
Thus $qBq=B^Rq$, and $q$ commutes with every element of $B^R$.
Thus the formula
\begin{equation}
\label{sigma}
\s(b)=bq
\end{equation}
defines a homomorphism $\s$ of $B^R$ onto the
$C^*$-subalgebra $B^Rq$ of $M(B\times_\a Q)$.
We will deduce from \propref{newrosenberg} below that
$\s$ is in fact an isomorphism.

Let $\b$ be the action of $Q/R$ on $B^R$ obtained from $\a$.
It is easy to see that the
maps $\s$ and $\tau$ from Equations~\ref{sigma} and \ref{tau}
combine to form a
covariant homomorphism
$(\s,\tau)$
of the action $(B^R,Q/R,\b)$, and that the
integrated form
\begin{equation}
\label{theta}
\t:=\s\times \tau:B^R\times_\b Q/R\to q(B\times_\a Q)q
\end{equation}
is surjective.

In the special case $R=Q$, the following is the main result of
\cite{rosenberg}:

\begin{prop}
\label{newrosenberg}
Let $(B,Q,\a)$ be an action,
let $R$ be a compact normal subgroup of $Q$,
let $(B^R,Q/R,\b)$ be the associated action,
and let $q=\Chi_R$.
Then
the map $\t:B^R\times_\b Q/R\to q(B\times_\a Q)q$ from \eqref{theta}
is an isomorphism.
\end{prop}

\begin{proof}
By the discussion preceding the statement of the proposition, it
remains to verify that $\t$ is injective,
and we do this by showing that for every covariant representation
$(\pi,U)$ of $(B^R,Q/R,\b)$ on a Hilbert space $V$ there exists a
representation $\rho$ of $q(B\times_\a Q)q$ on $V$ such that
$\rho\circ\t=\pi\times U$.

Recall from the theory of Rieffel induction \cite{rie:induced} that
the conditional expectation $\P:B\to B^R$ 
gives rise to a $B^R$-valued
inner product 
\[\<b,c\>_{B^R}=\P(b^*c)\]
on $B$,
so
the completion $X$ is
a Hilbert
$B^R$-module.
Moreover, $B$ acts on the left of $X$
by adjointable operators,
so we can use $X$ to induce $\pi$ to a representation
$\wilde\pi$ of $B$ on $\wilde V:=X\otimes_{B^R} V$.
An easy computation shows that the formula
\[
\wilde U_s(b\otimes \x)=\a_s(b)\otimes U_{sR}\x
\righttext{for}s\in Q,b\in B,\x\in V
\]
determines a representation $\wilde U$ of $Q$ on $\wilde V$ such
that $(\wilde\pi,\wilde U)$ is a covariant representation of
$(B,Q,\a)$.

Thus $\wilde\pi\times \wilde U$ is a representation of the
crossed product $B\times_\a Q$ on $\wilde V$; let $\rho_1$ be
its restriction to the corner $q(B\times_\a Q)q$.
We have $\rho_1(q)\wilde V=B^R\otimes_{B^R} V$,
because if $b\in B$ and $\x\in V$ then
\begin{align*}
\rho_1(q)(b\otimes \x)
&=\int_R\wilde U_r(b\otimes \x)\,dr
\\&=\int_R\bigl(\a_r(b)\otimes U_{rR}\x\bigr)\,dr
\\&=\int_R\a_r(b)\,dr\otimes \x.
\end{align*}
The subspace $B^R\otimes_{B^R}V$ is invariant for the
representation $\rho_1$;
let $\rho_2$ denote the associated subrepresentation of
$q(B\times_\a Q)q$.
A routine computation shows that
\[
W(b\otimes \x)=\pi(b)\x
\righttext{for}b\in B^R,\x\in V
\]
determines a unitary map $W$ of $B^R\otimes_{B^R}V$ onto $V$
which implements an equivalence between the representations
$\rho_2\circ \t$ and $\pi\times U$.
Thus we can take $\rho=\ad W\circ \rho_2$.
\end{proof}

\begin{cor}
\label{fix iso}
Let $(B,Q,\a)$ be an action,
let $R$ be a compact normal subgroup of $Q$,
let $(B^R,Q/R,\b)$ be the associated action,
and let $q=\Chi_R$.
Then the map $\s:B^R\to B^Rq$
from \eqref{sigma} is an isomorphism.
\end{cor}

\begin{proof}
It remains to observe that $\s$ is faithful, being
the composition of the injective
homomorphism $\t$ with the canonical embedding of $B^R$ into
$M(B^R\times_\b Q/R)$.
\end{proof}

\subsection*{Two projections}

If $A$ is a $C^*$-algebra and
$p$ is a projection in $M(A)$, then one of
the most basic applications of Rieffel's theory \cite{rie:induced} is
that 
the ideal $\bar{ApA}$
is Morita-Rieffel equivalent to 
the corner $pAp$ 
via the $\bar{ApA}-pAp$ imprimitivity bimodule $Ap$.
For
later purposes, we will need
a slightly more
subtle variant:

\begin{lem}
\label{morita}
Let $A$ be a $C^*$-algebra, and let $p,q\in M(A)$ be projections
with $p\le q$. Then
$q\bar{ApA}q$ is Morita-Rieffel equivalent to $pAp$.
\end{lem}

\begin{proof}
Just apply the above Morita-Rieffel equivalence $\bar{ApA}\sim pAp$
with $A$ replaced by $qAq$.
\end{proof}

\subsection*{Central projection}

Let $\b$ be an action of a locally compact group $T$ on a
$C^*$-algebra $C$, and let $d\in M(C)$ be a central projection. Then
$d$ may also be regarded as a multiplier of the crossed product
$C\times_\b T$, and it generates the ideal
\[\bar{(C\times_\b T)d(C\times_\b T)}.\]

\begin{prop}
With the above notation, we have:
\begin{enumerate}
\item $\bar{(C\times_\b T)d(C\times_\b T)}=I\times_\b T$,
where $I$ is the $T$-invariant ideal of $C$ generated by $d$.
\item $I=\clspn\{\b_t(d)C:t\in T\}=\{c\in C:p_\infty c=c\}$,
where $p_\infty=\sup\{\b_t(d):t\in T\}$.
\end{enumerate}
\end{prop}

\begin{proof}
(i)
follows from \cite{gre:local}*{Propositions~11~(ii) and 12~(i)}.

(ii)
The first equality
holds because $d$ is a central projection.
For the second, note that
the projections $\{\b_t(d):t\in T\}$ are central, so their supremum
$p_\infty$
is an open central projection in $C^{**}$, 
and the desired equality follows from, \eg,
\cite{ped}*{Proposition~3.11.9}.
To make this part of the
proof self-contained, we include the argument:
put
\[J=\{c\in C:p_\infty c=c\}.\]
For any $t\in T$ and $c\in C$ we have $\b_t(d)\le p_\infty$, so
\[p_\infty \b_t(d)c=\b_t(d)c.\]
Thus $I\subset J$. Suppose $a\in J$ but $a\notin I$. Then there
exists a 
nondegenerate
representation $\pi$ of $C$ such that $\pi(a)\ne 0$ but
$I\subset \ker\pi$. Extend $\pi$ to a weak*-weak-operator continuous
representation $\bar\pi$ of $C^{**}$. 
Enlarge the set $\{\b_t(d):t\in T\}$ to an
upward-directed set $P$ of central projections in $M(C)$, so that
there is an increasing net $\{p_i\}$ in $P$ converging weak* to
$p_\infty$.
Then
$p_ia\to p_\infty a$
weak*, so $\pi(p_ia)\to \bar\pi(p_\infty a)$. We have
$\bar\pi(p_\infty a)=\pi(a)$ because $a\in J$,
and $\pi(p_ia)=0$
for all $i$, so we deduce that $\pi(a)=0$,
a contradiction.
\end{proof}

\begin{q}
When will $p_\infty$ be a multiplier of $B^R$? (\exref{ex:Heisenberg} shows that it is not always so.)
\end{q}

\subsection*{Combined Results}

With the notation and assumptions of Proposition~\ref{newrosenberg}, 
put
\[A=B\times_\a Q.\]
Also let $d\in M(B)$ be an $R$-invariant central projection, so that
$d$ is also a central projection in $M(B^R)$.
Put
\[p_\infty=\sup\{\a_s(d):s\in Q\}.\]
Then $p_\infty$ is an open central projection in $(B^R)^{**}$.
Let $I$ be the $Q/R$-invariant ideal of $B^R$ generated by $d$.
We have $dq=qd\in M(A)$, and we denote this projection by $p$.

The following theorem combines the previous results in this section:

\begin{thm}
\label{combine}
With the above notation, we have:
\begin{enumerate}
\item $\t(I\times_\b Q/R)=q\bar{ApA}q$.
\item $I=\clspn\{\a_s(d)B^R:s\in Q\}=\{b\in B^R:p_\infty
b=b\}$.
\item $\s(I)=\clspn\{sds\inv qBq:s\in Q\}=\clspn\{sqdBqs\inv:s\in
Q\}$.
\item $pAp$ is Morita-Rieffel equivalent to $I\times_\b Q/R$.
\end{enumerate}
\end{thm}

\begin{proof}
The only part that still requires proof is (iii). We have
\begin{align*}
\s(I)
&=\clspn_{s\in Q}\t\circ\a_s(dB^R)
\end{align*}
because $\a_s(B^R)=B^R$.
For each $s\in Q$ we have
\begin{align*}
\s\circ\a_s(dB^R)
&=\s\circ\b_{sR}(dB^R)
\\&=\tau_{sR}\s(dB^R)\tau_{sR}^*
\righttext{(covariance)}
\\&=(sq)dB^Rq(sq)^*
=sqdB^Rqs\inv
=sqdqBqs\inv
\\&=sqdBqs\inv
\righttext{($dq=qd$)}
\\&=sds\inv qBq
\righttext{($sq=qs,sB=Bs$)},
\end{align*}
and (iii) follows.
\end{proof}

\section{Hecke crossed products}
\label{hecke crossed}

In this section our main object of study is a 
Schlichting pair $(G,H)$ 
which has 
the semidirect-product decomposition
of \eqref{setup}, with the additional condition that 
$R$ be normal in $Q$.
We shall obtain crossed-product
$C^*$-algebras which are Morita-Rieffel equivalent to the
completion of the Hecke algebra inside $C^*(G)$, 
similarly
to
certain results of \cite{hecke}.
At the end of the section we shall briefly indicate how our
results can be
applied if the Hecke pair is incomplete.

Put $A=C^*(G)$ and $B=C^*(N)$,
and let $\a$ denote the canonical action of $Q$ on $B$
determined by conjugation of $Q$ on $N$.
Then $A$ is isomorphic to the crossed product
$B\times_\a Q$, and
we identify these two $C^*$-algebras.

Normalize the Haar measures on $N$ and $Q$ so that $M$ and $R$ each
have measure $1$. Then the product measure is a Haar measure on $G$,
and
$H$ has measure $1$.
Thus $p_M:=\Chi_{M}$
is a central projection in $B$, hence is a projection in $M(A)$.
Similarly,
$p_R:=\Chi_{R}$
is a central projection in $C^*(Q)$, hence also a projection in
$M(A)$, and
we have
\[
p_H:=\Chi_{H}=p_Mp_R=p_Rp_M\in A.
\]

By \cite{hecke}*{Corollary~4.4}
the Hecke algebra of the pair $(G,H)$ is $\H=p_HC_c(G)p_H$, whose
closure in $A$ is the corner $p_HAp_H$.
From \secref{crossed} we get isomorphisms
\begin{align*}
\t&=\s\times\tau:B^R\times_\b Q/R\xrightarrow{\cong} p_RAp_R\\
\s&:B^R\xrightarrow{\cong} B^Rp_R\\
\tau&:C^*(Q/R)\xrightarrow{\cong} C^*(Q)p_R,
\end{align*}
and an ideal
\[I=\{b\in B^R:p_\infty b=b\}\ideal B^R,\]
where
\[p_\infty=\sup \{\a_s(p_M):s\in Q\}\in (B^R)^{**}.\]
\thmref{combine} quickly gives
the following analogue of
\cite{hecke}*{Theorem~8.2}:

\begin{thm}
\label{cross thm}
With the above notation:
\begin{flalign}
&\t(I\times_\b Q/R)=p_R\bar{Ap_HA}p_R.\tag{i}
\\&I=\clspn\{\a_s(p_M)B^R:s\in Q\}.\tag{ii}
\\&\s(I)
=\clspn\{sp_Ms\inv p_RBp_R:s\in Q\}\tag{iii}
\\&\qquad=\clspn\{sp_Rp_MBp_Rs\inv:s\in Q\}\notag
\\&\qquad=\clspn\{sp_Ms\inv p_Rnp_R:s\in Q,n\in N\}.\notag
\\&\text{$p_HAp_H$ is Morita-Rieffel equivalent to $I\times_\b
Q/R$.}\tag{iv}
\end{flalign}
\end{thm}

\begin{proof}
The only thing left to prove is the last equality of part (iii), and
this
follows from \thmref{combine}, because
$M$ is compact open in $N$, hence
\[p_MB=\clspn\{p_Mn:n\in N\}\]
(note that the projection $d$ from \thmref{combine} is $p_M$ here).
\end{proof}

\begin{rem}
Note that if $R$ is nontrivial then $p_H$ is never full in $A$:
Since $N$ is normal in $G$ with $Q=G/N$,
there is a natural homomorphism
$C^*(G)\to C^*(Q)$
which maps $p_H$ to $p_R$.
Thus $p_R$ is a nontrivial projection,
which, being central,
is not full in $C^*(Q)$.
\end{rem}

We say that the family $\{sMs\inv:s\in Q\}$ of conjugates of $M$ is
\emph{downward-directed} if the intersection of any two of them
contains a third.

\begin{prop}\label{down}
If $\{sMs\inv:s\in Q\}$ is downward-directed, then
\[p_R\bar{Ap_HA}p_R=p_RAp_R\cong B^R\times_\b Q/R.\]
\end{prop}

\begin{proof}
Because the pair $(G,H)$ is reduced we have
\[\bigcap_{s\in Q}sMs\inv=\{e\},\]
so the upward-directed set $\{sp_Ms\inv:s\in Q\}$ of projections has
supremum $p_\infty=1$ in $(B^R)^{**}$.
Therefore the ideal $I$ from \thmref{cross thm} coincides with $B^R$,
and the result follows.
\end{proof}

\begin{rem}
In the above proposition, we have
\[p_R\bar{Ap_HA}p_R=p_RAp_R\]
although the ideal $\bar{Ap_HA}$ of $A$ is
proper if $R$ is nontrivial.
\end{rem}

As in \cite{hecke}*{Section~7},
we specialize to the case where $N$ is abelian.
Taking Fourier transforms, the action $\a$ of $Q$ on $B$ becomes
an action $\a'$ on $C_0(\what N)$:
\[
\a'_s(f)(\p)=f(\p\circ \a_s)
\righttext{for}s\in Q,f\in C_0(\what N),\p\in \what N.
\]
The smallest $Q$-invariant subset of $\what N$ containing $M^\perp$
is
\[
\O=\bigcup_{s\in Q}(sMs\inv)^\perp.
\]
The Fourier transform of the fixed-point algebra $B^R$ is
isomorphic to
$C_0(\what N/R)$,
where
$\what N/R$
is the orbit space
under the action of $R$.
The smallest $Q/R$-invariant subset of $\what N/R$ containing
$M^\perp/R$ is $\O/R$.
Thus the Fourier transform of the ideal $I$ of $B^R$
is $C_0(\O/R)$.
Let $\g$ be the associated action of $Q/R$ on $C_0(\O/R)$.
The following corollary is analogous to
\cite{hecke}*{Corollary~7.1}.

\begin{cor}
With the assumptions and notation of Proposition~\ref{down},
if $N$ is abelian then $p_HAp_H$ is Morita-Rieffel equivalent to
the crossed product
$C_0(\O/R)\times_\g Q/R$.
\end{cor}

We finish this section with a brief indication of how the above
general theory can be used when $(G,H)$ is the Schlichting
completion of a reduced Hecke pair $(G_0,H_0)$.
More precisely, we assume that $G_0=N_0\rtimes Q_0$, $M_0\ideal
N_0$, $R_0\ideal Q_0$, $R_0$ normalizes $M_0$, and that
$(G_0,H_0)$ is a reduced Hecke pair (and
Propositions~\ref{hecke}--\ref{reduced}
give conditions under which the
latter happens). 
By
\corref{semidirect completion},
the closures
$N$, $Q$, $M$, and $R$ of 
$N_0$, $Q_0$, $M_0$, and $R_0$, respectively, satisfy the
conditions of the current section.
The action $(B,Q,\a)$ restricts to an action $(B,Q_0,\a_0)$, and
by density we have $B^R=B^{R_0}$.
The map $sR_0\mapsto sR$ for $s\in R_0$ gives an isomorphism
$Q_0/R_0\cong Q/R$ of discrete groups, and the action $\b$ of
$Q/R$ on $B^R$ corresponds to an action $\b_0$ of
$Q_0/R_0$ on $B^{R_0}$.
Thus we have a natural isomorphism
\[
B^R\times_\b Q/R\cong B^{R_0}\times_{\b_0} Q_0/R_0.
\]
Again by density,
for all $s\in Q$
there exists $s_0\in Q_0$ such that $p_Rs=p_Rs_0$,
and similarly for all $n\in N$
there exists $n_0\in N$ such that $np_M=n_0p_M$.
We deduce:

\begin{cor}
Using the above isomorphisms and identifications:
\begin{enumerate}
\item $I$ is the $Q_0/R_0$-invariant ideal of $B^{R_0}$
generated by $p_M$;
\item $I\times_{\b_0} Q_0/R_0\cong p_R\bar{Ap_HA}p_R$;
\item $p_\infty=\sup\{sp_Ms\inv:s\in Q_0\}$;
\item $I\cong \clspn\{sp_Ms\inv p_Rnp_R:s\in Q_0,n\in N_0\}$;
\item $p_HAp_H$ is Morita-Rieffel equivalent to $I\times_{\b_0}
Q_0/R_0$.
\end{enumerate}
\end{cor}

As explained in \cite{hecke}, many of the nice properties of the Hecke
algebra in \cite{bc} hold because the family $\{xHx\inv\mid
x\in G\}$ of conjugates of $H$ is downward directed; in particular this
implies that the projection $p$ is full. In our situation we can only
have $p$ full if $R=\{e\}$, but we do have the following:

\begin{cor}
\label{M directed}
Suppose the conjugates $\{sMs\inv\mid s\in Q\}$ of $M$ are downward
directed. Then  $I=B^{R_0}$ and $p_HAp_H$ is Morita-Rieffel equivalent
to $ B^{R_0}\times_{\b_0}Q_0/R_0$.
\end{cor}
\begin{proof}
We have $sp_Ms\inv = p_{ sMs\inv }$, so by the assumptions $p_\infty=1$.
\end {proof}

Continuing with $(G,H)$ being the Schlichting completion of
$(G_0,H_0)$ as above,
we again consider the special case where $N$, equivalently $N_0$,
is abelian.
Fourier transforming,
by density we have
\[
\O=\bigcup_{s\in Q_0}(sMs\inv)^\perp,
\]
and there is an associated action $\g_0$ of $Q_0/R_0$ on
$C_0(\O/R)$, giving:

\begin{cor}
With the above notation,
$p_HAp_H$ is Morita-Rieffel equivalent to
$C_0(\O/R)\times_{\g_0} Q_0/R_0$.
\end{cor}

\section{Examples}
\label{example}

We shall here illustrate the results from the preceding sections with
a number of examples.  
Some arguments are only sketched.

First note that the case $R=\{e\}$ is treated in \cite{hecke}*{Sections~7--8}.

\begin{ex}
\label{ex:general.bad} The situation with $M=\{e\}$ and $R\ideal Q$
is also interesting. From \secref{group} we see  that $(NQ,R)$ is
Hecke if and only if $R_{n,\{e\}}=\{r\in R\mid rnr\inv=n\}$ 
has finite index in $R$ for all $n$. 
The pair is reduced if and only if $\bigcap_n R_{n,\{e\}}=\{e\}$,
\ie, if the map $R\to \aut N$ is injective. Here $\NN=N$,
$p:=p_H=p_R$,  and \thmref{cross thm}
gives Morita-Rieffel
equivalences between $\bar{ApA}$, $pAp$ and $C^*(N)^R\times Q/R$.
\cite{hecke}*{Example~10.1} is a special case of this situation.
\end{ex}

We shall next study  $2\times 2$ matrix groups (and leave it to the
reader to see how this generalizes to $n\times n$ matrices).
For any ring $J$ we let $\M(2,J)$ denote the set of all $2\times 2$
matrices with entries in $J$; 
we let $\GL(2,J)$ denote the group of invertible
elements of $\M(2,J)$; 
$\SL_{\pm}(2,J)$ denotes the subgroup of
$\GL(2,J)$ consisting of those matrices with determinant $\pm 1$, and 
$\SL(2,J)$ is the subgroup of $\GL(2,J)$ of matrices with determinant $1$.

\begin{prop}
\label{N bar no p}
Suppose $N=\Q^2$, $M=\Z^2$,  $Q$ is a subgroup of 
$\GL(2,\Q)$ containing the diagonal subgroup
    $D=\{
(\begin{smallmatrix}
\lambda &0\\
0&\lambda
\end{smallmatrix})
\mid\, \lambda\in \Q^\times\}$,
and 
        $R=Q\cap \GLtwoZ$.
Then $(NQ,MR)$ is a reduced Hecke pair, and the Schlichting completion
is given by
\[
\NN={\cc A}_f^2 ,\quad\MM={\cc Z}^2,\quad \RR=\invlim R/R(s), \midtext{
and } \QQ=\bigcup_{q\in Q/R} q\RR,
\]
where
$\QQ$ has the topology from $\RR$, i.e.,
$q_i\rightarrow e$ if and only
if
$q_i\in \RR$ eventually and $q_i\rightarrow e$ in $ \RR$.
\end{prop}

\begin{proof}
Given
$q\in Q$ there is an integer matrix $k\in D$ such that $kq\inv$ is an
integer matrix. From this it follows that $kq\inv \Z^2\subset \Z^2$
and therefore $kMk\inv \subset qMq\inv$. This implies that the sets
$\{k\Z^2\}$ are downward-directed and form a base 
at $e$ for the Hecke topologies
of $M$ and $N$, by \propref{subbase}. We also note that $\bigcap_k
kMk\inv=\bigcap_k k\Z^2=\{e\}$, by \propref{reduced}.
Thus $\NN={\cc A}_f^2$ and $\MM={\cc Z}^2$,
with ${\cc A}_f$ the finite adeles  and 
${\cc Z}$ the integers in ${\cc A}_f$.

Next, if $n\in N$ there exists $s\in\Z$ such that $sn\in M$.
Take
 $n_1=(\begin{smallmatrix}
1/s\\
0
\end{smallmatrix})$
and
$n_2=(\begin{smallmatrix}
0\\
1/s
\end{smallmatrix})$. By definition 
$r\in R_{n,M}$ if and only if $(r-I)n\in \Z^2$.
One  checks that $R_{n_1,M}\cap R_{n_2,M}\subset R_{n,M}$ and that
    \[
    R_{n_1,M}\cap R_{n_2,M}=
    \left\{r\in R\mid\, r-I \in \M(2, s\Z)\right\}.
\]
Call this subgroup $R(s)$; it is clearly a normal subgroup of finite
index in $R$.

Suppose $q=(\begin{smallmatrix}
a&b\\
c&d
\end{smallmatrix})\in Q$, and
without loss of generality we may assume  $q\in \M(2,\Z)$. 
Putting $t=\det (q)=ad-bc$, for $r\in R(t)$ we have
\[
q\inv( r- I)q= t\inv (\begin{smallmatrix}
d&-b\\
-c&a
\end{smallmatrix})
(r-I) (\begin{smallmatrix}
a&b\\
c&d
\end{smallmatrix})
\in \M(2,\Z),
\]
and it follows that $q\inv rq\in \M(2,\Z)$. The same argument holds for
$r\inv$, so
both $q\inv rq$ and $q\inv r\inv q$ are integer matrices in $Q$.
Thus
\[q\inv rq\in Q\cap \GLtwoZ=R.\]
From this it follows that
\[
R(t)\subset R\cap qRq\inv \quad\text{ for }\quad t=\det (q),
\]
and we have just observed that $[R:R(t)]<\infty$,
so $[R:R_q] < \infty$.

The same argument also shows that
$R(st)\subset R\cap qR(s)q\inv$ for any~$s$, and
therefore to any given finite sets $E\subset Q$ and $F\subset N$
there exists $s\in\N$ such that 
$R(s)\subset R^E_F$.
Combining all this with
\propref{subbase} we see that the family $\{R(s)\mid\, s\in \N\}$ is a
base
at $e$ for the Hecke topology restricted to  $R$ or~$Q$.

Finally, note that $\bigcap_s R(s)=\{e\}$.
\end{proof}

A similar result holds when $\Q$ is replaced by other number fields,
\eg, $\Z[p\inv]$ for a prime number $p$
(not to be confused with the projection~$p$).
We state it without proof:

\begin{prop}
\label{N bar with p}
Suppose $N=\Z[p\inv]^2$, $M=\Z^2$,    $Q$ is a subgroup of $\GL
(2,\Z[p\inv])$ containing the diagonal subgroup
$D=\{
(\begin{smallmatrix}
p^n&0\\
0&p^n
\end{smallmatrix})
\mid\, n\in \Z\}$,
and 
        $R=Q\cap \GLtwoZ$.
 Then $(NQ,MR)$ is a reduced Hecke pair, and the Schlichting completion
is given by
\[
\NN=\Q_p^2, \quad\MM={\cc Z}_p^2,\quad
\RR=\invlim R/R(p^n),\midtext{ and }\QQ=\bigcup_{q\in Q/R}q\RR,
\]
where as above
$\QQ$ has the topology from $\RR$.
\end{prop}

\begin{ex}
\label{ex:GL2} 

Let us first consider the maximal $p$-adic  case with $Q=\GL
(2,\Z[p\inv])$ and $R= \GLtwoZ$.

\begin{prop} Let
\label{G=TH-p-adic}
$T=\{
(\begin{smallmatrix}
p^m&0\\
c&p^n
\end{smallmatrix})
\mid\, m,n\in\Z,c\in \Z[p\inv]\}$.
Then
$T\,\SL_\pm (2,{\cc Z}_p)=
\{
g\in \GL (2,\Q_p)\mid\, \det(g)\in \pm p^\Z \}$.
\end{prop}
\begin{proof}
Clearly the left hand side is included in the right hand side.
For the opposite it suffices to show that every
$g\in \M(2,{\cc Z}_p)$ with $\det(g)\in p^\N$
is a member of the left hand side.
Let $g=(\begin{smallmatrix}a&b\\c&d\end{smallmatrix})$.

Case 1. Suppose
$b=0$ and
$ad=p^m$. If $a=p^nu$ with $u$ a unit in ${\cc Z}_p$, we must
have
$d=u\inv p^{m-n}$. So
    \[
g=\begin{pmatrix}
a&0\\
c&d
\end{pmatrix}
=\begin{pmatrix}
p^n&0\\
0&p^{m-n}
\end{pmatrix}
\begin{pmatrix}
1&0\\
x&1
\end{pmatrix}
\begin{pmatrix}
u&0\\
0&u\inv
\end{pmatrix}
\]
with $x=cu\inv p^{n-m}$.
Now $x=y+z$ with $y\in \Z[1/p]$ and 
$z\in {\cc Z}_p$,
and since
$(\begin{smallmatrix}
1&0\\
z&1
\end{smallmatrix})
(\begin{smallmatrix}
u&0\\
0&u\inv
\end{smallmatrix})  \in \SL (2,{\cc Z}_p)$ it follows that
$g=(\begin{smallmatrix}
a&0\\
c&d
\end{smallmatrix})\in T\,\SL (2,{\cc Z}_p)$.

Case 2. Suppose $a=0$ and $b\neq 0$. Then
    \[
    g=\begin{pmatrix}
0&b\\
c&d
\end{pmatrix}
=\begin{pmatrix}
b&0\\
d&-c
\end{pmatrix}
\begin{pmatrix}
0&1\\
-1&0
\end{pmatrix}
\in T\SL (2,{\cc Z}_p).
\]

Case 3. Suppose
$a=p^mu$ and $b=p^nv$ with $u,v$ units in ${\cc Z}_p$. We may assume
$m\geq n$, if not we multiply 
by
$(\begin{smallmatrix}
0&1\\
-1&0
\end{smallmatrix})$ as in Case 2.
So $p^{-n}a\in {\cc Z}_p$. Then
    \[
\begin{pmatrix}
a&b\\
c&d
\end{pmatrix}
=\begin{pmatrix}
p^{n}&0\\
v\inv d&p^{-n}ad-vc
\end{pmatrix}
\begin{pmatrix}
p^{-n}a&v\\
-v\inv&0
\end{pmatrix}.
\]
The second matrix on the right hand side is in $\SL (2,{\cc Z}_p)$,
while the first has determinant equal to $ad-bc$ which by assumption
is in $p^\N$,
so by Case~1 
this matrix is in $T\,\SL (2,{\cc Z}_p)$.
\end{proof}

\begin{thm}
\label{GL2p} Let $Q=\GL (2,\Z[p\inv])$ and $R= \GLtwoZ$.
Then
\begin{enumerate}
\item $\RR=\invlim R/R(p^n)=\SL _{\pm}(2,{\cc Z}_p)$;
\item $\QQ=\bigcup_{q\in Q/R}q\SL _{\pm}(2,{\cc Z}_p)=
\{
g\in \GL (2,\Q_p)\mid\, \det(g)\in \pm p^\Z \}$,
\end{enumerate}
where $\QQ$ has the topology from
$\RR= \SL _{\pm}(2,{\cc Z}_p)$.
\end{thm}

\begin{proof}
Since $Q=TR$ we get $\QQ=T\RR$, which by \propref{G=TH-p-adic} equals
the
right hand side.
That 
\cite{hecke}*{Theorem~3.8} applies to the pair
$(\NN\QQ,\MM\RR)$
now follows from 
Propositions~\ref{N bar with p} and
\ref{G=TH-p-adic},
and density of
$\GLtwoZ$ in
$\SL _{\pm}(2,{\cc Z}_p)$
(see \cite{kri}*{Proposition~IV.6.3}).
\end{proof}

Now look at the case 
where
$Q=\GL (2,\Q)$ and $R= \GLtwoZ$. 
We first need a version of \propref{G=TH-p-adic}:

\begin{prop} Let
\label{G=TH}
$T=\{
(\begin{smallmatrix}
a&0\\
c&d
\end{smallmatrix})
\mid\, a,c,d\in\Q, ad\neq 0\}$.
Then $T\,\SL (2,{\cc Z})=
\left\{
g\in \GL (2,{\cc A}_f)\mid\, \det(g)\in \Q\right\}$.
\end{prop}
\begin{proof}
Again one inclusion is obvious, so suppose 
$g=(\begin{smallmatrix}
a&b\\
c&d
\end{smallmatrix})\in \GL (2,{\cc A}_f)$ with $\det g\in\Q$, 
in fact without loss of generality
we may
assume $\det g=1$.
For each prime $p$ let 
$g_p=(\begin{smallmatrix}
a_p&b_p\\
c_p&d_p
\end{smallmatrix})$ be the corresponding matrix in $\GL (2,\Q_p)$. 
For all but finitely many $p$ we will have 
$g_p\in \SL (2,{\cc Z}_p)$. In these cases take $k_p=g_p$.

In the other cases we can not have both $a_p$ and $b_p$ zero, 
so by \propref{G=TH-p-adic} there is a matrix 
$k_p\in \SL (2,{\cc Z}_p)$ such that 
$g_pk_p\inv\in T\cap \GL (2,\ \Z[1/p])$.
So $k=(k_p)\in \SL (2,{\cc Z})$ and $gk\inv\in T$ as claimed.
\end{proof}

\begin{thm}
\label{GL2}
Let $Q=\GL (2,\Q)$ and $R=\GLtwoZ$.
Then
\begin{enumerate}
\item $\RR=\SL _{\pm}(2,{\cc Z})$;
\item $\QQ=\bigcup_{q\in Q/R}q\SL _{\pm}(2,{\cc Z})=
\left\{
g\in \GL (2,{\cc A}_f)\mid\, \det(g)\in \Q\right\}$,
\end{enumerate}
where $\QQ$ has the topology from
$\RR= \SL _{\pm}(2,{\cc Z})$.
\end{thm}

\begin{proof}
From \cite{kri}*{Proposition~IV.6.3} (the hard part is hidden there)
it follows that
\begin{align*}
\RR
=\invlim R/R(s)
=\invlim\SL_{\pm}(2,\Z_s)
=\SL _{\pm}(2,{\cc Z}).
\end{align*}

Since $\QQ=T\RR$, (ii) follows from \propref{G=TH}.
\end{proof}

Note that the topology on $\QQ$ is not the relative topology from $\GL
(2,{\cc A}_f)$, in contrast with \thmref{GL2p}.

This is 
essentially
the same result as
\cite{lln:hecke}*{Proposition~2.5}.
Since $R$ is not normal in $Q$ we can not use
\thmref{cross thm}, but it would be interesting
to get a description of the $C^*$-algebra $p_RAp_HAp_R$ in these
cases (see \cite{connes-marcolli}). However, note that we are not
using exactly the
same algebra, since in both \cite{connes-marcolli} and
\cite{lln:hecke} the action of $Q$ is by left multiplication on
$\M(2,\Q)$.

\end{ex}

\begin{ex}
\label{ex:ax+b}
Much recent work on Hecke algebras started with the study of the
affine group over $\Q$ in \cite{bc}.
Other number fields have also been extensively studied,
as in,
\eg, \cite{connes-marcolli} and
\cite{laca-franken}.
For a survey, see \cite{connes-marcolli}*{Section~1.4}.
We shall here illustrate how our approach works for a quadratic
extension of $\Q$. For details about the number theory used here we
refer to the book \cite{nzm}.

Let $d$ be a square-free integer 
such that%
\footnote{If, for instance, $d=5$ one
should instead use  $M=\Z[(1+\sqrt 5)/2]$,
\etc\ (see \cite{nzm}*{Theorem~9.20}).}
$d\not\equiv 1 \mod 4$,
and let $N=\Q(\sqrt d)$,
$M=\Z[\sqrt d]$, $Q=\Q(\sqrt d)^\times$,
and $R=\{r\in Q\,\mid\,r,r\inv\in M\}$.

So
\[
    R=\{m+n\sqrt d\,\mid\,m,n\in \Z, m^2-dn^2=\pm 1\}
\]
is the group of units in the field $N$.
An alternative matrix description is as follows:
\begin{align*}
N&=\Q^2, \qquad M=\Z^2, \\
Q&=\left\{
\begin{pmatrix}
a&db\\
b&a
\end{pmatrix}
\biggm|
a,b\in \Q, a^2-db^2\neq 0\right\},\\
R&=\left\{
\begin{pmatrix}
m&dn\\
n&m
\end{pmatrix}
\biggm|
m,n\in \Z, m^2-dn^2=\pm 1\right\}.
\end{align*}
So we  get $\NN={\cc A}_f^2$ and $\MM={\cc Z}^2$. 

Here \thmref{cross thm} applies, so
    \[
    pAp\sim_{MR} C_0({\cc A}_f^2/\RR)\rtimes Q/R.
\]
In this way we obtain
\cite{laca-franken}*{Proposition~3.2}
for  the field  $\Q(\sqrt d)$
without using the theory of semigroup crossed products, and this will
also work in greater generality.

The structure of these crossed products can be studied by the
Mackey-Takesaki orbit method as in
\cite{lr:ideal};
note that the
orbit closures in $ \NN/\RR$ under the action of $ Q/R$
are basically the same as the orbit closures in $\NN$ under the action
of $ Q$.

To determine $\RR$ and its topology we need some more information.
First, if $d<0$ then $R$ is finite (of order 2 or 4).
So let us concentrate on the case with $d>1$.
Then, by
\cite{nzm}*{Theorem~7.26}
we have
$R\cong \{\pm 1\}\times\Z$, and 
in fact there exists $r_0\in R$ such that $R=\{\pm
r_0^n\,\mid\,n\in\Z\}$.
For instance if $d=2$ one can take $r_0=1+\sqrt 2$.

Let us 
look at $R(s)$.
There is a smallest integer $n_s>0$ such that $r_0^{n_s}\equiv 1 \bmod
s$.
From this we get
$\RR=\invlim R/R(s)=\{\pm 1\}\times \invlim \Z/\Z_{n_s}$.
However, examples show that the behavior of the numbers $n_s$ is
complicated, so a more exact description of $\RR$ is difficult.

Perhaps counterintuitively, in general it turns out that
\[
 \RR\subsetneq\{ m+n \sqrt d
\mid\,m,n\in {\cc Z}, m^2-dn^2=\pm 1\}.
\]
This is because under the homomorphism $\Z[\sqrt d]\mapsto \Z_s[\sqrt
d]$ the units $R$ in $\Z[\sqrt d]$ are in general mapped onto a proper
subgroup of the units in $\Z_s[\sqrt d]$. For instance $4$ is a unit in 
$\Z_{17}[\sqrt 2]$, but 
$\pm (1+\sqrt 2)^n\not\equiv 4 \mod 17$ for all $n$.
\end{ex}

\begin{ex}
\label{ex:Heisenberg}
We shall here give a slightly different treatment of the Heisenberg
group
than in \cite{hecke}. Take
\begin{align*}
N&=\Q/\Z\times\Q, 
&M&=\{0\}\times\Z,\\
Q&=\left\{
\begin{pmatrix}
1&q\\
0&1
\end{pmatrix}
\biggm| q\in \Q\right\},
&R&=\left\{
\begin{pmatrix}
1&r\\
0&1
\end{pmatrix}
\biggm| r\in \Z\right\},
\end{align*}
with the obvious action of $Q$ on $N$. If
$x=
(\begin{smallmatrix}
1&1/n\\
0&1
\end{smallmatrix})$ with $n\in \N$ one checks that
$M\cap xMx\inv=\{0\}\times n\Z$.
So we have
\[
\NN=\Q/\Z\times {\cc A}_f={\cc A}_f/{\cc Z} \times {\cc A}_f
\midtext{ and }
\MM=\{0\}\times{\cc Z}.
\]
If
$n=(\begin{smallmatrix}
a\\
b/m
\end{smallmatrix})$
with $b,m\in\Z$ and
$r=
(\begin{smallmatrix}
1&r\\
0&1
\end{smallmatrix})$,
then $rnr\inv -n\in M$ if and only if $rb\in m\Z$. Thus
\[
\QQ=\left\{
\begin{pmatrix}
1&q\\
0&1
\end{pmatrix}
\biggm| q\in {\cc A}_f\right\}
\midtext{ and }
\RR=\left\{
\begin{pmatrix}
1&r\\
0&1
\end{pmatrix}
\biggm| r\in {\cc Z}\right\}.
\]
We have $\widehat\NN={\cc Z}\times {\cc A}_f$
and $\MM^\perp={\cc Z}\times{\cc Z}$. Moreover, the dual action of $\QQ$
on $\widehat\NN$ is given by
\[
(z,w)
\begin{pmatrix}
1&q\\
0&1
\end{pmatrix}= (z,qz+w).
\]
\begin{lem}
\begin{align*}
\Omega:&=\bigcup_{q\in \QQ} q\MM^\perp
=\left\{(z,qz+w)\,\mid\, z,w\in {\cc Z}, q\in {\cc A}_f\right\}\\
&=\left\{(z,u) \in {\cc Z}\times {\cc A}_f\,\mid\, z_p=0 \implies
u_p\in {\cc Z}_p \right\}.
\end{align*}
\end{lem}
\begin{proof}
Clearly if $(z,w)\in \Omega$ and $z_p=0$, then $w_p\in {\cc Z}_p$.

Conversely, suppose $(z,u)$ is an element of the right hand side. If
$u_p\in {\cc Z}_p $, take $q_p=1$ and $w_p=u_p-z_p\in {\cc Z}_p$.
For the finitely many $p$ with  $u_p\notin {\cc Z}_p $, we have
$u_p=x_p+v_p$ with $x_p\in\Q^\times$ and $v_p\in {\cc Z}_p $, and by
assumption $z_p\neq 0$. Take $q_p =z_p\inv x_p\in\Q_p$, so $q_p
z_p+w_p=u_p$. Thus with $q:=(q_p)\in{\cc A}_f$ and $w:=(w_p)\in{\cc
Z}$, we have  $q z+ w= u$.
\end{proof}
So here $\Omega$ is open but not closed, hence the projection $p_\infty$
defined in \secref{hecke crossed} is not in $M(B^R)$.

The orbits under the action of $R$ can be described  as follows: $(0,w)$
is always a fixed point.
If $z\neq 0$, then the $R$-orbit of $(z,w)$
is
$(z, w+z{\cc Z})$.
\end{ex}

\end{document}